\documentclass[journal]{IEEEtran}
\usepackage{algorithm,algcompatible,algpseudocode}
\usepackage{amsmath}

\usepackage{amssymb}
\usepackage{dirtytalk}
\usepackage{enumitem}
\setlist[enumerate]{leftmargin=.5in}
\setlist[itemize]{leftmargin=.5in}

\usepackage{amsopn}

\newcommand{\I}{\mathcal{I}}

\newcommand{\E}{\mathbb{E}}
\newcommand{\D}{\mathcal{D}}
\newcommand{\yk}{y^k_{N_k}}
\newcommand{\ve}{\varepsilon}
\newcommand{\e}{\epsilon}
\newcommand{\be}{\beta}

\let \oldsection \section
\renewcommand{\section}{\vspace{3ex plus 1ex}\oldsection}
\newcommand{\BEAS}{\begin{eqnarray*}}
	\newcommand{\EEAS}{\end{eqnarray*}}
\newcommand{\BEA}{\begin{eqnarray}}
\newcommand{\EEA}{\end{eqnarray}}

\newcommand{\BEQ}{\begin{equation}}
\newcommand{\EEQ}{\end{equation}}
\newcommand{\BIT}{\begin{itemize}}
	\newcommand{\EIT}{\end{itemize}}
\newcommand{\BNUM}{\begin{enumerate}}
	\newcommand{\ENUM}{\end{enumerate}}

\newcommand{\BA}{\begin{array}}
	\newcommand{\EA}{\end{array}}

\newcounter{dummy} \numberwithin{dummy}{section}
\newtheorem{mythm}[dummy]{Theorem}

\newtheorem{mylem}[dummy]{Lemma}

\newtheorem{myh}{A.}

\newcommand{\mr}{\mathrm}
\newcommand{\mb}{\mathbb}
\newcommand{\mc}{\mathcal}
  \usepackage{subfig}

\usepackage{cite}
\ifCLASSINFOpdf
   \usepackage[pdftex]{graphicx}
\else
\fi
\begin{document}
\title{A Fast Stochastic Plug-and-Play ADMM for Imaging Inverse Problems}

\author{Junqi Tang,
       % Karen Egiazarian,~\IEEEmembership{Fellow,~IEEE}
        %and~
        Mike~Davies,~\IEEEmembership{Fellow,~IEEE}% <-this % stops a space
\thanks{J.Tang and M.Davies are with the School of Engineering, University of Edinburgh. Correspondence to J.Tang@ed.ac.uk}% <-this % stops a space
}


\markboth{June, 2020}%
{Shell \MakeLowercase{\textit{et al.}}: Bare Demo of IEEEtran.cls for IEEE Journals}



\maketitle

\begin{abstract}
 In this work we propose an efficient stochastic plug-and-play (PnP) algorithm for imaging inverse problems. The PnP stochastic gradient descent methods have been recently proposed and shown improved performance in some imaging applications over standard deterministic PnP methods. However, current stochastic PnP methods need to frequently compute the image denoisers which can be computationally expensive. To overcome this limitation, we propose a new stochastic PnP-ADMM method which is based on introducing stochastic gradient descent inner-loops within an inexact ADMM framework. We provide the theoretical guarantee on the fixed-point convergence for our algorithm under standard assumptions. Our numerical results demonstrate the effectiveness of our approach compared with state-of-the-art PnP methods.
\end{abstract}


% \fi
\IEEEpeerreviewmaketitle

\section{Introduction}

Recent trends in the research of computational imaging have been focusing on developing algorithms which are able to jointly utilize the power of classical physical models and advanced image priors \cite{venkatakrishnan2013plug, egiazarian2007compressed,kamilov2017plug, romano2017little,reehorst2018regularization}. These methods typically take the form of well-known optimization algorithms, and plug in a pretrained deep neural network \cite{zhang2017beyond} or a patch-based denoiser with non-local denoising properties \cite{buades2005non,dabov2008image,chen2017trainable,talebi2013global}. In this work we propose a novel stochastic plug-and-play method for imaging inverse problems. Consider the following observation model for a linear inverse problem:
\begin{equation}\label{eq:1}
b = A x^\dagger + w, \ \ \ A \in \mb{R}^{n \times d}
\end{equation}
where $x^\dagger$ denotes the vectorized (raster) ground truth image, $A \in \mb{R}^{n \times d}$ represents the forward measurement model, $b \in \mb{R}^n$ denotes the observation, while $w \in \mathbb{R}^n$ represents the random additive noise. Traditionally, in order to get a good estimate of the ground truth $x^\dagger$, we typically seek to find the minimizer of a composite objective function:
\begin{equation}\label{reg_ls}
    x^\star \in \arg \min_{x \in \mathbb{R}^d} \left\{ f(x) + g(x)\right\},
\end{equation}
where $f(x) = \frac{1}{n}\sum_{i = 1}^n f_i(x) = \frac{1}{n}\sum_{i = 1}^n f(a_i, b_i, x)$ is the data fidelity term which is assumed to be convex and smooth, such as the least-square loss $ f(x) = \frac{1}{2n}\|Ax - b\|_2^2$, and here we denote $a_i$ as the i-th row of $A$ and $b_i$ as the i-th element of $b$. Meanwhile $g(x)$ in (\ref{reg_ls}) denotes a regularization term which encodes image priors, with classical examples including the sparsity-inducing regularization in wavelet domain and the total-variation regularization \cite{chambolle2016introduction}, etc.  The composite loss function (\ref{reg_ls}) can be effectively minimized via a class of iterative algorithms which are known as the proximal splitting methods \cite{combettes2011proximal}, including the forward-backward splitting \cite{lions1979splitting,2009_Beck_Fast,beck2009fast}, primal-dual splitting\cite{chambolle2011first,chambolle2015convergence,pesquet2014class}  and the Douglas-Rachfort splitting/alternating direction method of multipliers (ADMM)\cite{douglas1956numerical, boct2015inertial, boyd2011distributed}, etc.

Considering the link between proximal operators and denoising, researchers \cite{venkatakrishnan2013plug,egiazarian2007compressed,kamilov2017plug} have discovered that if they simply replace the proximal operator on $g$ with a direct call to an off-the-shelf denoising algorithm, such as NLM \cite{buades2005non}, TNRD \cite{chen2017trainable}, BM3D \cite{dabov2008image}, or the DnCNN \cite{zhang2017beyond}, excellent image recovery/reconstruction results can often be attained. Although the imaging community has very limited theoretical understanding and convergence analysis for such algorithms so far, these ad-hoc approaches have shown state-of-the-art performances in various imaging applications.  %, with small demands on the training data compared to pure deep neural-network approaches which are based on end-to-end training.
 
 Recently, inspired by the success of stochastic gradient descent (SGD) methods in solving large-scale optimization tasks in machine learning \cite{bottou2010large,shalev2011pegasos,kingma2014adam} and some imaging applications \cite{chambolle2018stochastic,chouzenoux2017stochastic,ehrhardt2017faster}, Sun et al \cite{sun2019online} have extended the deterministic plug-and-play ISTA/FISTA method \cite{kamilov2017plug} and proposed PnP-SGD in order to improve the computational efficiency. In each iteration of PnP-SGD, a minibatch stochastic gradient is computed as an unbiased estimator of the full gradient, which yields a computational benefit in each iteration. However, as discussed in \cite{tang2019practicality}, current stochastic gradient methods in general need to compute the proximal operator/denoisers more frequently than the deterministic gradient methods within the same amount of gradient evaluations. When the denoiser is computationally expensive, the actual performance benefit of using stochastic gradient techniques may be compromised due to this computational overhead. In this work, we seek to address this issue for stochastic PnP methods and propose a more practical approach for utilizing the power of stochastic gradient techniques to accelerate the deterministic plug-and-play algorithms.

 \subsection{Main contributions}
 This paper's contribution is two-fold:
 
 \begin{itemize}
     \item  We propose an efficient stochastic PnP method which empirically improves upon previous stochastic approach in \cite{sun2019online} by reducing the number of calls on the modern denoisers which are usually a bottleneck for computation. We demonstrate the effectiveness of our approach in X-ray computed tomography (CT) imaging problems.
     \item We provide a theoretical fixed-point convergence analysis for our stochastic PnP algorithm under standard assumptions.
 \end{itemize}

 \section{Stochastic PnP-ADMM}

In a recent work of Sun et al \cite{sun2019online}, a stochastic PnP algorithm is proposed for image restoration and reconstruction, which can be written as the following:
  \begin{eqnarray*}
 && \mathrm{\textbf{PnP-SGD}} - \mathrm{Initialize}\ x^0, z^0 = x^0 \\
 &&\mathrm{For} \ \ \ k = 1, 2,...,  t\\
&&\left\lfloor
\begin{array}{l}
x^{k}=  \mathcal{D} [z^{k-1}-\eta \cdot \triangledown f_{S_k}(z^{k-1})]\\
z^{k} = x^{k} + \alpha^k(x^k - x^{k-1}) 
\end{array}
\right.
 \end{eqnarray*}
 where $\triangledown f_{S_k}(x^{k-1})$ denotes a minibatch stochastic gradient with a randomly subsampled index $S_k$, chosen uniformly at random from a partitioned index $\hat{\I} = \{\I_1, \I_2, ..., \I_K\}$, where $\I_1\cup\I_2\cup...\cup\I_K = [n]$ and $\I_i \cap \I_j = \emptyset$, $\forall i \neq j \in [K]$. The PnP-SGD algorithm is essentially a plug-and-play variant of the stochastic proximal gradient descent \cite{rosasco2014convergence} which is based on the forward-backward splitting \cite{lions1979splitting}. In each iteration of PnP-SGD, a minibatch stochastic gradient estimate is computed:
\begin{equation}\label{sg}
    \triangledown f_{S_k}(z^k) = \frac{1}{m}\sum_{i \in S_k}  \triangledown f_i(z^k),
\end{equation}
where $m = \frac{n}{K}$. It first performs a stochastic gradient descent step with a step-size $\eta$, then a denoising step is computed using an off-the-shelf denoiser \cite{buades2005non,dabov2008image,chen2017trainable,zhang2017beyond} denoted as $\D(\cdot)$. Finally a momentum step is performed for empirical convergence acceleration with a momentum parameter $\alpha_k$ as in the FISTA algorithm \cite{2009_Beck_Fast}. The computational benefit of PnP-SGD over its deterministic counterparts (PnP-ISTA/PnP-FISTA \cite{kamilov2017plug}) comes from using an approximation of the full gradient by the minibatch gradient (\ref{sg}) which can be efficiently computed. However, the PnP-SGD does not have the capability to reduce the cost of computing the denoising step -- the denoiser has to be called at each iteration. To overcome this computational bottleneck, one plausible approach is to decouple the gradient step and the denoising step via Douglas-Rachford splitting/ADMM instead of the forward-backward splitting. In this work we study and propose a stochastic gradient extension of the PnP-ADMM algorithm\footnote{We write the update rule of PnP-ADMM in this paper using the equivalent Douglas-Rachford splitting reformulation \cite[Section 9.1]{ryu2019plug} for the simplicity of notation in analysis.} \cite{venkatakrishnan2013plug}:
 \begin{eqnarray*}
 && \mathrm{\textbf{PnP-ADMM}} - \mathrm{Initialize}\ x^0=z^0 \in \mb{R}^d; \\
 &&\mathrm{For} \ \ \ k = 0, 2,...,  t\\
&&\left\lfloor
\begin{array}{l}
y^{k}=  \mathrm{prox}_{\tau f} [z^k]\\
x^{k+1} = \mathcal{D}[2 y^k - z^{k}]\\
z^{k+1} = z^{k} + x^{k+1} - y^k.
\end{array}
\right.
 \end{eqnarray*}
where in each iteration an exact proximal step on the data-fidelity term $f(x)$ is computed with a constant step-size $\tau$:
  \begin{equation}\label{prox_step}
    \mathrm{prox}_{\tau f} (\cdot) = \arg \min_{x} \frac{1}{2} \|x - \cdot \|_2^2 + \tau f(x).
\end{equation}
In classical ADMM the step-size $\tau$ can be any positive constant to ensure convergence. The update rule of the PnP-ADMM can be written as $z^{k+1} = T(z^k)$, where $T(\cdot)$ is an operator defined as \cite{ryu2019plug}:
\begin{equation}
    T = \frac{1}{2}I + \frac{1}{2}(2\mc{D} - I)(2\mr{prox}_{\tau f} - I).
\end{equation}
Our proposed solution, presented in algorithm 1, is to use SGD with momentum to \emph{approximately} solve the prox step (\ref{prox_step}) within PnP-ADMM framework. We denote the number of inner-iterations at the $k$-th outerloop as $N_k$. In each inner-iteration a stochastic gradient descent step is performed with a step-size $\eta_k$, and then followed by a momentum step for empirical acceleration. Unlike the PnP-SGD which needs to call the denoiser in every stochastic gradient descent iteration, the proposed method only needs to compute the denoiser once every $N_k$ iterations. Our theoretical analysis is restricted to the case where we set $\alpha_j = 0$ and the parameters $N_k$ and $\eta_k$ are chosen adaptively in each outer-iteration. However, in practice, a constant step size $\eta$ which is inversely proportional to the Lipschitz constant $\eta_k = O(\frac{1}{L})$, $\tau = O(1)$, a constant number of inner-iterations $N_k = O(K)$, and a FISTA-like momentum parameter $\alpha_j = \frac{j-1}{j+3}$ \cite{chambolle2015convergence} are suggested for good empirical performance.

 \begin{algorithm}[t]\label{AA}
   \caption{--- Stochastic PnP-ADMM}
\begin{algorithmic}

   \State  {\bfseries Initialization:} number of inner-iterations: $[N_1, N_2, ...,N_K]$, momentum parameter sequence: $[\alpha_1, \alpha_2,...,\alpha_{\max_{j\in[K]} N_j}]$, partition index $\hat{\I} = \{\I_1,\I_2,...,\I_K\}$, $z^0 \in \mb{R}^d$, $y_0^0 \in \mb{R}^d$.
   \For{$k = 1$ {\bfseries to} $K$}
   \For{$j =1$ {\bfseries to} $N_k$}
    \State Randomly sample $S_j \in \hat{\I}$ with replacement.
    \State Compute a stochastic gradient estimator $\triangledown f_{S_j}(y^k_{j-1})$
   \State        $v^k_{j}=  y^k_{j-1}-\eta_k \cdot [\tau \triangledown f_{S_j}(y^k_{j-1}) + y^k_{j-1} - z^k]$

    \State {\bfseries Momentum:} $y^k_{j} = v^k_{j} + \alpha_j (v^k_{j} - v^k_{j-1})$
   
   \EndFor
   \State $x^{k+1} = \mathcal{D}( 2 y^k_{N_k} - z^k)$
   \State $z^k = z^{k-1} + x^{k+1} - y^{k}_{N_k}$
   \State $y^{k+1}_0 = x^{k+1}$
   \EndFor
\State Output $x^{K}$
\end{algorithmic}
\end{algorithm} 

\section{Convergence Analysis}

In this section we provide theoretical analysis for our stochastic PnP-ADMM. When we run a stochastic gradient-based innerloop, we are effectively making an approximation of the proximal step $y^{k}=  \mathrm{prox}_{\tau f} [z^{k}]$, hence we can write our stochastic PnP-ADMM algorithm as the inexact recursion:
\begin{equation}\label{DR}
    z^{k+1} = T(z^k) + \ve^k.
\end{equation}
where $\ve^k$ denotes the approximation error. Now the desired fixed-point convergence can be established for (\ref{DR}), if we make the following standard assumptions as in \cite{ryu2019plug} on the denoiser and the data-fidelity term.

\subsection{Generic Assumptions}

\begin{myh}
 The denoiser satisfies:
 \begin{equation}
     \|(\D - I)(x) - (\D - I)(y)\|_2 \leq \be \|x - y\|_2, \ \forall x, y \in \mb{R}^d,
 \end{equation}
 with $\be > 0$.
\end{myh}
It is easy to show that A.1 implies a relaxed non-expansiveness condition on $\D(\cdot)$ which reads $\|\D(x) - \D(y)\|_2 \leq (1+ \be) \|x - y\|_2$ and is satisfied for a wide class of modern denoisers such as NLM and properly trained DnCNNs \cite{ryu2019plug}. 

\begin{myh}
 $f(\cdot)$ is $\mu$-strongly-convex:
 \begin{equation}
     f(x) - f(y) - \langle \triangledown f(y), x-y \rangle \geq \mu \|x - y\|_2^2,\ \forall x, y \in \mb{R}^d,
 \end{equation}
 with $\mu > 0$. Meanwhile for a given minibatch partition index $\hat{\I} = \{\I_1, \I_2, ..., \I_K\}$ such that $f(\cdot) = \frac{1}{K}\sum_{k=1}^K f_{\I_k}(\cdot)$, each $f_{\I_k}$ is $L$-smooth, such that $\forall x, y \in \mb{R}^d$:
  \begin{equation}
     f_{\I_k}(x) - f_{\I_k}(y) - \langle \triangledown f_{\I_k}(y), x-y \rangle \leq L \|x - y\|_2^2.
 \end{equation}
\end{myh}

The strong-convexity assumption is necessary for our analysis. It seems pessimistic since a number of imaging inverse problems do not have strong-convexity. We believe that the assumption on strong-convexity could be relaxed, e.g. following the ideas from \cite{oymak2017sharp,bolte2007lojasiewicz,liang2017local,pmlr-v70-tang17a,tang2018rest}. On the other hand, one can instead run Algorithm 1 on a regularized objective $\hat{f}(x) = f(x) + \frac{\e}{2}\|x\|_2^2$ to manually enforce strong-convexity, which is a classical trick in convex optimization \cite{nesterov2013introductory}. Nevertheless, we believe that relaxing this assumption is an important future direction for the analysis.

\subsection{Analysis}

We first apply an existing convergence result for SGD for establishing the approximation accuracy of the inner-loop:
\begin{mylem}\label{lem_sg}
Under Assumption A.2, denote that for $k$-th outer-loop of Stochastic PnP-ADMM $y_\star^{k}=  \mathrm{prox}_{\tau f} [z^k]$, and define the following quantities for each outer-iteration $k$:
\begin{equation}
\begin{aligned}
   & \sigma_k^2 := \E_q[\|\tau\triangledown f_{\I_q}(y^k_\star) + y^k_\star - z^k\|_2^2], \\& \xi_k := \|y_\star^{k} - x^{k}\|_2^2,
    \end{aligned}
\end{equation}
then if the step size $\eta_k = \frac{\mu_0 \varepsilon}{2 \varepsilon \mu_0 L_0 + 2\sigma_k^2}$, $\alpha_j = 0$ for all $j$, $N_k = 2\log(\frac{\xi_k}{\ve})(\frac{L_0}{\mu_0} + \frac{\sigma_k^2}{\mu_0^2 \ve})$ with $\mu_0 = \tau \mu + 1$, $L_0 = \tau L + 1$, then we have the approximation error of the proximal step bounded as:
\begin{equation}
    \E \|y_{N_k}^k - \mathrm{prox}_{\tau f} [z^k]\|_2^2 \leq \varepsilon,
\end{equation}
where the expectation is taken over the random sampling of the indices within the inner-loop.
\end{mylem}

{\bf Proof.} 
We first observe that the proximal step $y_\star^{k}=  \mathrm{prox}_{\tau f} [z^k]$ can be written precisely as a finite-sum optimization problem of the follow form:
\begin{equation}
\begin{aligned}
    \mr{prox}_{\tau f}(z^k) &= \arg \min_x \frac{1}{2}\|x - z^k\|_2^2 + \tau f(x)\\
    &= \arg \min_x \frac{1}{K}\sum_{q=1}^K [\tau f_{\I_q}(x) + \frac{1}{2}\|x - z^k\|_2^2],
\end{aligned}
\end{equation}
which is a $(\tau \mu + 1)$-strongly-convex objective and each of the element in the sum is $(\tau L + 1)$-smooth.

According to \cite[Theorem 2.1]{needell2014stochastic}, if we run SGD (starting at $x^k \in \mb{R}^d$) with uniform random sampling and a  step size $\frac{\mu_0\ve}{2\ve\mu_0 L_0 + 2\sigma_k^2}$, then after $N = 2\log(\frac{\|y_\star^k - x^k\|_2}{\ve})(\frac{L_0}{\mu_0} + \frac{\sigma_k^2}{\mu_0^2 \ve})$ iterations, we have $\E \|y_{N_k}^k - y_\star^k\|_2^2 \leq \varepsilon$. $\ \ \ \ \ \ \ \ \ \ \ \ \ \ \ \ \ \ \ \ \ \ \ \Box$

Now, we are able to prove the fixed-point convergence for the inexact recursion (\ref{DR}), and hence for our proposed method.

\begin{mythm}\label{main_thm2}
Assume A.1 and A.2 with $\be < 1$, denote positive values $\sigma$ and $\xi$ such that $\forall k$,  $\sigma_k \leq \sigma$ and $\xi_k \leq \xi$, and the quantities $\mu_0$, $L_0$, $\sigma_k^2$, $\xi_k$ are defined as in Lemma \ref{lem_sg}. If we choose the step-size parameters as $\tau > 1/(1+\be-2\be^2)$, $\alpha_j = 0$, $\eta_k = \frac{\mu_0}{2 \mu_0 L_0 + 2k\sigma^2}$, $N_k = 2\log(k\xi)(\frac{L_0}{\mu_0} + \frac{k\sigma^2}{\mu_0^2})$, we have the following fixed-point convergence for Algorithm 1:
\begin{equation}
    \E\|z^{k+1} - z^k\|_2 \rightarrow 0,
\end{equation}
when $k \rightarrow +\infty$.
\end{mythm}

Our main theorem suggests that for the basic form of Algorithm 1 where we choose the momentum $\alpha_j = 0$, the outerloop step-size $\tau = O(1)$, the inner-loop step-size decreasing in $k$ and the number of inner-iterations increasing in $k$, Algorithm 1 is guaranteed to converge to a fix point. However our numerical results in section IV suggest that we may set the number of inner-loop $N_k$ and step size $\eta_k$ to be constant and use FISTA-type of momentum \cite{2009_Beck_Fast,beck2009fast,chambolle2015convergence} for good empirical performance in practice.

\subsection{Proof for Theorem \ref{main_thm2}}
Firstly, due to assumption A.1 we have:
\begin{equation}
    \begin{aligned}
     \be \|x - y\|_2 &\geq \|(\D - I)(x) - (\D - I)(y)\|_2 \\&\geq \|\D(x) - \D(y)\|_2 - \|x - y\|_2,\ \forall x, y \in \mb{R}^d,
     \end{aligned}
\end{equation}
hence $\|\D(x) - \D(y)\|_2 \leq (1+ \be) \|x - y\|_2 ,\ \forall x, y \in \mb{R}^d$. Denote $u_k = \mr{prox}_{\tau f} (z^k) - \yk = y_\star^k - \yk$, we have:
\begin{equation}
    \begin{aligned}
        \|\ve^k\|_2 &= \|z^{k+1} - T(z^k)\|_2\\
                    &\leq \|u_k\|_2 + \|\D(2\yk - z^k) - \D(2\yk - z^k + 2 u_k) \|_2\\
                    &\leq (3+2\be)\|u_k\|_2.
    \end{aligned}
\end{equation}
Applying Lemma \ref{lem_sg} gives $ \E \|\ve^k\|_2 \leq \frac{3+2\be}{k}$. Now according to \cite[Theorem 2]{ryu2019plug}, under assumption A.1 and A.2, we can ensure that:
\begin{equation}
    \|T(x) - T(y)\|_2 \leq \delta \|x - y\|_2,\ \forall x, y \in \mb{R}^d,
\end{equation}
where $\delta = \frac{1 + \be + \be\tau\mu + 2 \be^2\tau\mu}{1 + \tau\mu +2\be\tau\mu}$. Moreover, if $\tau > 1/(1+\be-2\be^2)$ and $\be < 1$, then $\delta < 1$. Hence we have:
\begin{equation}
\begin{aligned}
     \E\|z^{k+1} - z^k\|_2 &= \E\|T(z^{k}) - T(z^{k-1}) + \ve^k - \ve^{k-1}\|_2 \\
     &\leq \E\|T(z^{k}) - T(z^{k-1})\|_2 \\
      &\ \ +\E\|\ve^k\|_2 +\E\| \ve^{k-1}\|_2\\
           &\leq \delta \E\|z^k - z^{k-1}\|_2 + \frac{3 + 2\beta}{k} + \frac{3+2\beta}{k-1}\\
     &\leq \delta \E\|z^k - z^{k-1}\|_2 + \frac{15}{k}
\end{aligned}
\end{equation}
If we recursively apply the same argument we will get:
\begin{equation}
    \E\|z^{k+1} - z^k\|_2 \leq \delta^k \E\|z^1 - z^0\|_2 + \frac{15}{k}\sum_{i = 0}^{k-1} \frac{k\delta^i}{k - i}.
\end{equation}
Then we use a classic criterion to show the boundedness of series $\sum_{i=0}^{+\infty} v_i$ where $v_i = \frac{k\delta^i}{k - i}$. For any finite $p > 1$, we have:
\begin{equation}
    \lim_{i \rightarrow + \infty} i^pv_i = \lim_{i \rightarrow + \infty} \frac{\delta^i i^p k}{k - i} = 0,
\end{equation}
where we take $k \rightarrow + \infty$ and $i \leq k-1$, and then:
\begin{equation}
    \sum_{i = 0}^{k-1} \frac{k\delta^i}{k - i} < + \infty, \ \ \ \ \frac{15}{k}\sum_{i = 0}^{k-1} \frac{k\delta^i}{k - i} \rightarrow 0.
\end{equation}
Hence by taking $k \rightarrow +\infty$, we have $\E \|z^{k+1} - z^k\|_2 \rightarrow 0$. Thus finishes the proof for Theorem \ref{main_thm2}.

\begin{figure}[t] % %{\textwidth}
	\centering	
    \includegraphics[width=.469\textwidth]{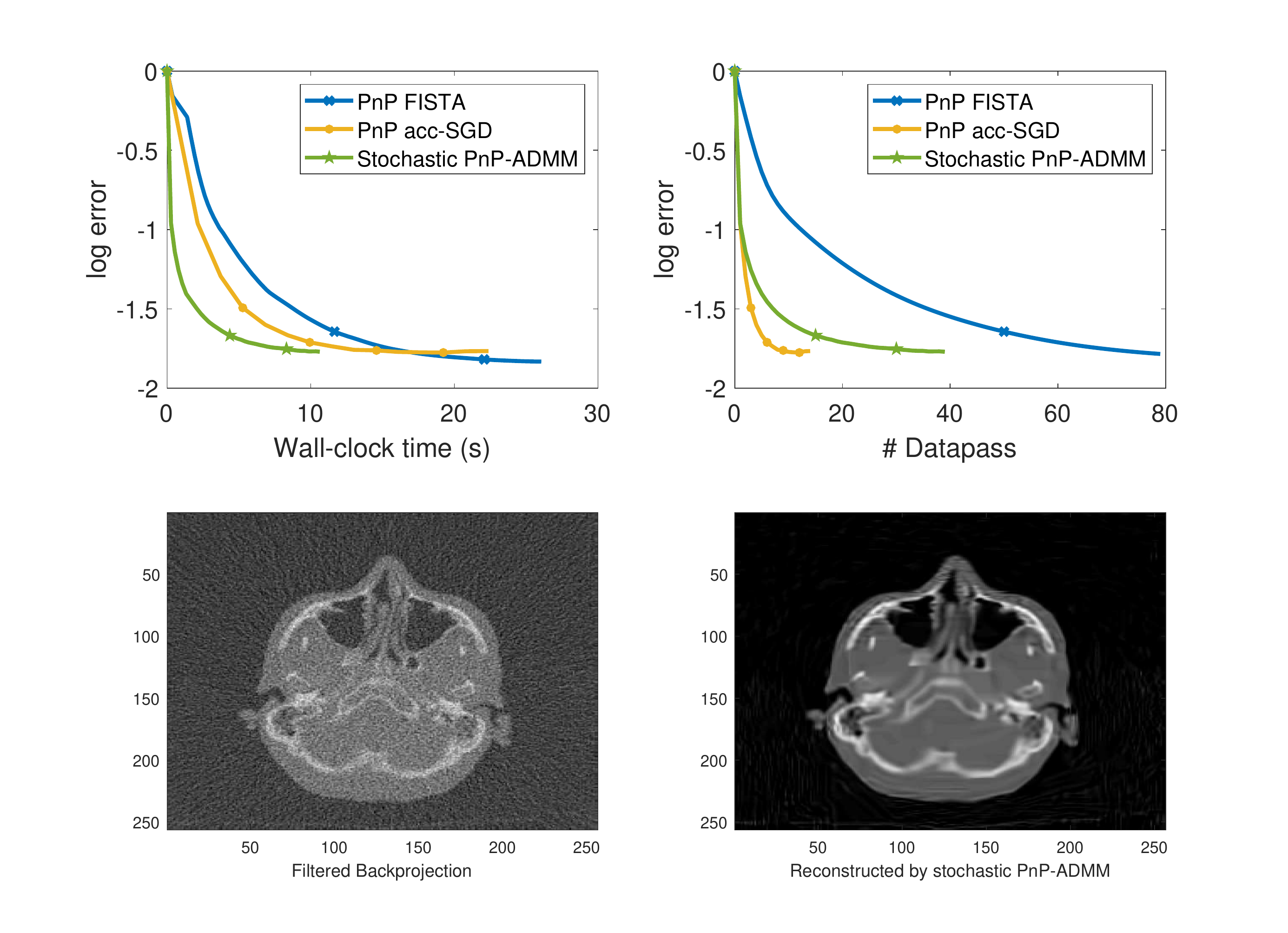}
	\caption{The estimation error plot of the compared algorithms on low-dose CT example}
	\label{fig:1}
\end{figure}

\section{Numerical Experiments}

 For our numerical experiments, we choose the X-ray CT imaging as an example since it is know to favor stochastic gradient methods \cite{tang2019practicality}. We compare our algorithm with the state-of-the-art stochastic PnP method with momentum acceleration proposed by Sun et al \cite{sun2019online}, as well as the PnP FISTA algorithm \cite{kamilov2017plug}. We use MATLAB R2018a in a machine with 1.6 GB RAM, 1.80 GHz Intel Core i7-8550 CPU. 
 
 We first test the compared methods on low-dose CT imaging problems, where low-energy noisy X-ray measurements with $I_0 = 10^3$ are used, which demands strong image priors are used in order to achieve good-quality reconstructions. Meanwhile we also compare these algorithms in sparse-view CT imaging with $I_0 = 10^{4}$, where fewer X-ray measurements are taken compared to the number of pixels to be inferenced\footnote{The numerical result in this example suggests that empirically the strong-convexity is not needed for the stochastic PnP-ADMM to converge.}. For low-dose CT example, we choose the penalized weighted least-squares objective as the data-fidelity term, which is tailored for low-dose CT \cite{wang2006penalized}. For sparse-view CT example, we choose the standard least-squares loss as the data fidelity term. For the randomized methods we partition the data into 10 minibatches. The noisy CT observations are obtained via $y \sim \mathrm{Poisson}(I_0 e^{- Ax})$ where the forward operator $A$ is implemented using the {\it AIRtools} package \cite{hansen2012air}. For our algorithm, we set $N_j = 10$ for all $j$ such that in each inner-loop we make exactly one pass of the data, outer-loop step-size $\tau = 1$, inner-loop step-size $\eta_k = \frac{1}{L}$, and the momentum parameter $\alpha_j = \frac{j-1}{j+3}$ as suggested in \cite{chambolle2015convergence}.
 
 \begin{figure}[t] % %{\textwidth}
	\centering	
    \includegraphics[width=.469\textwidth]{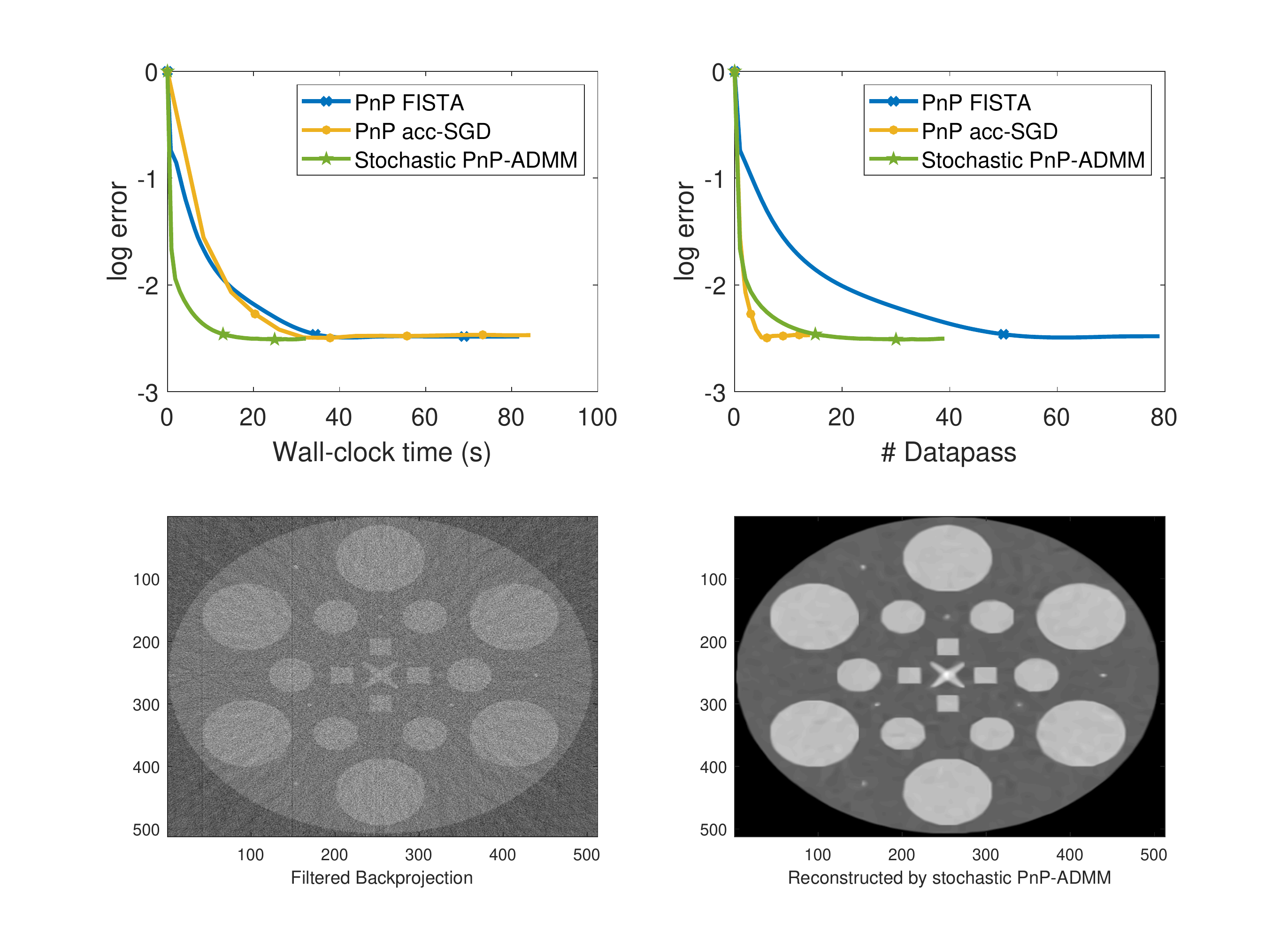}
	\caption{The estimation error plot of the compared algorithms on sparse-view CT example}
	\label{fig:3}
\end{figure}
 
We choose the BM3D \cite{dabov2008image} with the denoiser-scaling \cite{xu2020boosting} as the denoiser:
\begin{equation}
    \D_\gamma (x) = \frac{1}{\gamma} \mr{BM3D}(\gamma x),
\end{equation}
and maximize the reconstruction performance for each of the compared algorithms via grid-searching the parameter $\gamma$.

We present the numerical results of the algorithms in Figure 1 for low-dose CT inverse problem of size $A \in \mb{R}^{88256 \times 65536}$, and in Figure 2 for sparse-view CT imaging task of size $A \in \mb{R}^{92160 \times 262144}$. We plot the estimation error $\log_{10}\|x - x^\dagger\|_2$ to the ground-truth image $x^\dagger$, against the actual run time as well as the number of datapasses. We can observe that both PnP-SGD and our method are much faster than PnP-FISTA in terms of datapass. The PnP-SGD appears to be faster than our method in terms of number of datapasses. However, in terms of actual run time, the PnP-SGD is slower than our method, due to the need to compute the costly BM3D at each stochastic gradient iteration.

\section{Conclusion}

In this work we propose a stochastic PnP-ADMM algorithm which is able to provide practical acceleration with stochastic gradient techniques, for efficiently solving imaging inverse problems. This is an effective approach to make the stochastic PnP schemes truly practical, by reducing the computational overhead of the modern denoisers. We provide a fixed-point convergence analysis, and demonstrate the effectiveness of our method in numerical experiments.

\section*{Acknowledgment}
This work is supported by ERC Advanced grant 694888, C-SENSE.

\ifCLASSOPTIONcaptionsoff
  \newpage
\fi

\bibliographystyle{IEEEtran}
\bibliography{main.bib}

\end{document}